%% file: agt-1-37.tex
\documentclass{gtart}
\input agtout

\lognumber{37}
\volumenumber{1}
\volumeyear{2001}
\papernumber{37}
\pagenumbers{743}{762}
\received{11 October 2000}
\revised{17 July 2001}
\accepted{29 Novemver 2001}
\published{4 December 2001}

\usepackage{amsmath}
\usepackage{amsfonts}
\usepackage{amssymb}

\newtheorem{theorem}{Theorem}[section]

\newtheorem{corollary}[theorem]{Corollary}

\newtheorem{lemma}[theorem]{Lemma}
\newtheorem{technical lemma}[theorem]{Technical Lemma}

\newtheorem{proposition}[theorem]{Proposition}

\theoremstyle{remark}

\newtheorem{remark}[theorem]{Remark}

\begin{document}

\title{Splitting of Gysin extensions}
\author{A. J. Berrick\\A. A. Davydov}

\asciiaddress{Department of Mathematics\\
National University of Singapore
\\2 Science Drive 2, Singapore 117543, 
SINGAPORE\\and\\Department of Mathematics
\\Macquarie University, Sydney, NSW 2109
\\AUSTRALIA}

\address{Department of Mathematics\\
National University of Singapore
\\2 Science Drive 2, Singapore 117543, SINGAPORE}

\secondaddress{Department of Mathematics
\\Macquarie University, Sydney, NSW 2109
\\AUSTRALIA}

\email{berrick@math.nus.edu.sg, davydov@ics.mq.edu.au}

\begin{abstract}
Let $X\rightarrow B$ be an orientable sphere bundle. Its Gysin sequence
exhibits $H^{\ast}(X)$ as an extension of $H^{\ast}(B)$-modules. We prove that
the class of this extension is the image of a canonical class that we define
in the Hochschild $3$-cohomology of $H^{\ast}(B),$ corresponding to a
component of its $A_{\infty}$-structure, and generalizing the Massey triple
product. We identify two cases where this class vanishes, so that the Gysin
extension is split. The first, with rational coefficients, is that where $B$
is a formal space; the second, with integer coefficients, is where $B$ is a torus.
\end{abstract}

\asciiabstract{Let X --> B be an orientable sphere bundle. Its Gysin
sequence exhibits H^*(X) as an extension of H^*(B)-modules. We prove
that the class of this extension is the image of a canonical class
that we define in the Hochschild s-cohomology of H^*(B), corresponding
to a component of its A_infty-structure, and generalizing the Massey
triple product. We identify two cases where this class vanishes, so
that the Gysin extension is split. The first, with rational
coefficients, is that where B is a formal space; the second, with
integer coefficients, is where B is a torus.}

\asciikeywords{Gysin sequence, Hochschild homology, differential graded
algebra, 
formal space, A_infty-structure, Massey triple product}
\maketitle

\section{Introduction}

The paper is devoted to the description of the module structure of the
cohomology $H^{\ast}(X)$ of an oriented sphere bundle $X\rightarrow B$ over
the cohomology of the base $H^{\ast}(B)$.

This is a problem with a long history. In the late 'fifties (a decade or so
after the discovery of the Gysin sequence), several attempts were made to
determine multiplicative properties of the cohomology $H^{\ast}(X),$
principally by means of triple Massey products of $H^{\ast}(B)$ \cite{hirsch,
mas}. An obstacle to a full answer was the partial definiteness of Massey
products. W. S. Massey, in his review of a 1959 paper of D. G. Malm
\cite{malm}, summarized the situation as follows:

\textit{The author points out that it is misleading to say that `the
cohomology ring of the total space depends on the cohomology of the base space
and certain characteristic classes.' Just exactly what aspects of the homotopy
type of the base space it does depend on seems like a very complicated
question. Certainly various higher order cohomology operations are needed.}

Ironically, a mathematical framework for fully defined higher products
($A_{\infty}$-algebras) appeared just a little bit later, in a different
context, in the paper of J. Stasheff \cite{stash}. It was however only in 1980
that T. Kadeishvili defined an $A_{\infty}$-structure of a special kind for
the cohomology of a topological space or, more generally, of a differential
graded algebra \cite{kad}. Informally speaking, an $A_{\infty}$-structure on
cohomology is a collection of multilinear maps, the first of which is
multiplication satisfying certain compatibility conditions (see Remark
\ref{ashc} for the definition). The Massey products fit into this structure as
specializations of the multilinear maps. For example, the triple Massey
product is a specialization for Massey triples of the second component of the
$A_{\infty}$-structure (see Remark \ref{mass}). Another problem of higher
multiplications, their non-canonical nature, remained unhandled. For Massey
products this appears in the peculiar shape of the range of values; while for
$A_{\infty}$-structures it forces them to be defined only up to an isomorphism.

Here we attack the problem by using relations between $A_{\infty}$-structures
and Hochschild cohomology. The compatibility conditions imply that the second
component of an $A_{\infty}$-structure\ is a Hochschild $3$-cocycle, and the
corresponding Hochschild cohomology class is invariant under those
isomorphisms between $A_{\infty}$-structures that preserve ordinary
multiplication (see Remark \ref{ashc}). We give a self-contained construction
of this Hochschild $3$-cocycle for the cohomology of a differential graded
algebra, and prove the canonical nature of the corresponding Hochschild
cohomology class (Section \ref{sm}). As in \cite{kad}, we need some sort of
freeness condition on $H^{\ast}(B)=H^{\ast}(B;\,k)$ as a module over its
coefficient ring $k$. Throughout we assume that $k$ is commutative and
hereditary (for example, a principal ideal domain), so that every submodule of
a free module is projective.

The Gysin sequence provides an extension of $H^{\ast}(B)$-modules
\[
0\rightarrow H^{\ast}(B)/(cH^{\ast}(B)[\left|  c\right|  ])\rightarrow
H^{\ast}(X)\rightarrow\operatorname*{Ann}\nolimits_{H^{\ast}(B)}(c)[\left|
c\right|  -1]\rightarrow0,
\]
where $c\in H^{\ast}(B),$ of degree $\left|  c\right|  ,$ is the
characteristic class of the fibre bundle, and the square brackets indicate
grading shifts. We prove that the class of this extension in the extension
group
\[
\operatorname*{Ext}\nolimits_{H^{\ast}(B)}^{1}(\operatorname*{Ann}%
\nolimits_{H^{\ast}(B)}(c)[\left|  c\right|  -1],\,H^{\ast}(B)/(cH^{\ast
}(B)[\left|  c\right|  ]))
\]
is the image of our Hochschild $3$-cohomology class of the algebra $H^{\ast
}(B),$ under a natural homomorphism
\begin{align*}
HH_{\operatorname*{gr}}^{3}(H^{\ast}(B),\,\overline{H^{\ast}(B)}
[1])&\rightarrow\\
\operatorname*{Ext}\nolimits_{H^{\ast}(B)}^{1}&
(\operatorname*{Ann}\nolimits_{H^{\ast}(B)}(c)[\left|  c\right|
-1],\,H^{\ast}(B)/(cH^{\ast}(B)[\left|  c\right|  ])).
\end{align*}
To prove this, we use the language of differential graded (dg-) algebras. We
use the explicit construction of the Hochschild $3$-cohomology class
$[\theta]$ of the algebra $H^{\ast}(C)$ to interpret the class of the
$H^{\ast}(C)$-module extension
\[
0\rightarrow H^{\ast}(C)/(cH^{\ast}(C)[\left|  c\right|  ])\rightarrow
H^{\ast}(\operatorname*{cone}(l_{s(c)}))\rightarrow\operatorname*{Ann}%
\nolimits_{H^{\ast}(B)}(c)[\left|  c\right|  -1]\rightarrow0,
\]
obtained from the long exact sequence of the mapping cone
$\operatorname*{cone}(l_{s(c)})$ of left multiplication $l_{s(c)}:C[\left|
c\right|  ]\rightarrow C$, whenever $s(c)$ denotes a representative cocycle of
a cohomology class $c\in H^{\ast}(C)$ (section \ref{smm}). By expressing
Thom's Gysin sequence as the long exact sequence of the mapping cone
\cite{malm}, we are able to derive our main result (Theorem \ref{mco}).

\begin{theorem}
Let $S^{m}\rightarrow X\rightarrow B$ be an orientable sphere bundle, where
the homology $H_{\ast}(B;\,k)$ of the base is a projective $k$-module. Let
$c\in H^{m+1}(B)$ be its characteristic class and
\[
\eta\in\operatorname*{Ext}\nolimits_{H^{\ast}(B)}^{1}(\operatorname*{Ann}%
{}_{H^{\ast}(B)}(c)[m],\,H^{\ast}(B)/(cH^{\ast}(B)[m+1]))
\]
the class of the Gysin extension
\[
0\rightarrow H^{\ast}(B)/(cH^{\ast}(B)[m+1])\rightarrow H^{\ast}%
(X)\rightarrow\operatorname*{Ann}{}_{H^{\ast}(B)}(c)[m]\rightarrow0.
\]
Then the above natural homomorphism sends the Hochschild cohomology class of
$\theta(B)$ to the class $\eta$.
\end{theorem}

Applications of this theorem depend on calculation of the class $[\theta(B)]$
in specific cases. In the case of rational coefficients, there is a large and
well-studied family of spaces for which the obstruction is readily seen to
vanish, namely formal spaces. This gives the following theorem.

\begin{theorem}
Let $S^{m}\rightarrow X\rightarrow B$ be an orientable sphere bundle having
base $B$ a formal space, and characteristic class $c\in H^{m+1}(B)$. Then the
Gysin sequence with rational coefficients
\[
0\rightarrow H^{\ast}(B)/(cH^{\ast}(B)[m+1])\rightarrow H^{\ast}%
(X)\rightarrow\operatorname*{Ann}\nolimits_{H^{\ast}(B)}(c)[m]\rightarrow0
\]
splits as a sequence of right $H^{\ast}(B)$-modules.
\end{theorem}

On the other hand, the case of integer cohomology is far more delicate. There
seems to be no existing literature to which one can appeal. We attack the case
of a torus base, which is of fundamental importance to understanding the
multiplicative cohomology of integer Heisenberg groups, whose additive
structure has recently been determined \cite{lp}. In Section \ref{torus} below
we outline our proof of triviality of the canonical Hochschild $3$-cohomology
class for the torus, leading to the second application.

\begin{theorem}
Let $S^{m}\rightarrow X\rightarrow T$ be an orientable sphere bundle with
torus base $T$ and characteristic class $c\in H^{m+1}(T)$. Then the Gysin
sequence with integer coefficients
\[
0\rightarrow H^{\ast}(T)/(cH^{\ast}(T)[m+1])\rightarrow H^{\ast}%
(X)\rightarrow\operatorname*{Ann}\nolimits_{H^{\ast}(T)}(c)[m]\rightarrow0
\]
splits as a sequence of $H^{\ast}(T)$-modules.
\end{theorem}

Recently (Gda\'{n}sk conference, June 2001), D. Benson announced joint work
with H. Krause and S. Schwede concerning an obstruction, in $3$-dimensional
Hochschild cohomology, to decomposability of the cohomology of a dga-module.
Presumably their class will be closely related to ours; unfortunately, their
arguments seem to be no less lengthy. In applications to the Tate cohomology
of certain finite groups, their obstruction turns out to be nontrivial. That
finding is consistent with determination of nontrivial higher Massey products
with finite coefficients in \cite{Kraines}.

The results presented here were announced at the Arolla conference in August
1999, and the second-named author would like to thank D. Arlettaz for
providing the opportunity to speak there. Unfortunately, another journal's
prolonged inability to find a referee has led to a delay in publication. On
the other hand, we are grateful to this journal's referee for several
pertinent suggestions that have, we trust, made the paper more readable. In
particular, we have replaced the original six page technical proof of
triviality of the obstruction with torus base, by a brief summary of the
argument. The second author acknowledges the support of National University of
Singapore grant RP3970657.

\section{Secondary multiplication in cohomology}

\subsection{Hochschild cohomology of algebras}

By the \emph{Hochschild cohomology} of a unital, associative $k$-algebra $R,$
projective over a commutative ring $k,$ with coefficients in an $R$-bimodule
$M,$ we mean the extension groups
\[
HH^{\ast}(R,M)=\operatorname*{Ext}\nolimits_{R\text{-}R}^{\ast}(R,M)
\]
between $R$-bimodules $R$ and $M$. It can be identified \cite{CartanEilenberg}
(IX.6) with the cohomology of the complex
\begin{align*}
CH^{\ast}(R,M) =\operatorname*{Hom}\nolimits_{k}(R^{\otimes l},M),
\end{align*}
with differential
\[
\delta:CH^{l}(R,M)\rightarrow CH^{l+1}(R,M),
\]
\begin{align*}
(\delta a)(x_{1},...,x_{l+1})&= x_{1}a(x_{2},...,x_{l+1})+\\&
\sum_{i=1}^{l}(-1)^{i}a(x_{1},...,x_{i}%
x_{i+1},...,x_{l+1})+(-1)^{l+1}a(x_{1},...,x_{l})x_{l+1}.
\end{align*}

While it is well-known in this case that $\delta$ is indeed a differential,
for future reference we record the simple fact that establishes this property.
Its proof is a straightforward induction argument.

\begin{lemma}[The Differential Lemma]In the free associative ring
\[
\mathbb{Z}\left\langle u_{1},v_{1},u_{2},v_{2},\ldots\right\rangle ,
\]
for all $n,$%
\[
\left(  \sum_{i=1}^{n+1}u_{i}\right)  \left(  \sum_{i=1}^{n}v_{i}\right)
=\sum\limits_{1\leq i\leq j\leq n}(u_{i}v_{j}+u_{j+1}v_{i}).
\]
\end{lemma}

For the application at hand, it is easy to see that, whenever $i\leq j,$ we
have $\delta^{i}\delta^{j}=\delta^{j+1}\delta^{i},$ so we map $u_{i}$ and
$v_{i}$ to $(-1)^{i}\delta^{i}$ in order to obtain $\delta\delta=0.$

\remark
\label{ashc}$A_{\infty}$-structures and Hochschild cohomology
\cite{kad1}

\noindent The appearance of Hochschild cohomology in this paper is explained
by its connection with cohomological higher multiplications, namely with
$A_{\infty}$-algebras \cite{stash}. Recall that a minimal $A_{\infty}$-algebra
is a $k$-module $R=\bigoplus\nolimits_{s\geq0}R_{s}$ graded by nonnegative
integers, together with a collection of graded maps $m_{i}:R^{\otimes
i+2}\rightarrow R[i]$ for $i=0,1,...,$ satisfying
\begin{equation}
\sum_{i+j=n}\sum_{h=1}^{i+2}(-1)^{h(j-1)+(n+1)j}m_{i}\circ_{h}m_{j}=0,
\label{aia}%
\end{equation}
for any $n$ and suitably defined operations $\circ_{h}$(see also\cite{ger}).

Note that, for $n=0,$ the condition (\ref{aia}) implies that $(R,m_{0})$ is an
associative $k$-algebra. With $\overline{\delta}$ as the differential of the
standard Hochschild complex of the algebra $R=(R,m_{0})$ with coefficients in
the $R$-bimodule $\overline{R}$, the condition (\ref{aia}) for $n=1$ has the
form $\overline{\delta}(m_{1})=0$. This means that $m_{1}$ is a Hochschild
$3$-cocycle of $R$ with coefficients in $\overline{R}$. The cohomology class
$[m_{1}]\in HH^{3}(R,\overline{R})$ will be called the \emph{canonical class}
(or secondary multiplication) of the $A_{\infty}$-algebra $R$.

In the paper \cite{kad}, T Kadeishvili extended the ordinary multiplication on
the cohomology of a differential graded algebra to the minimal $A_{\infty}
$-algebra structure. In this paper we use only the second component of this
structure. For the sake of completeness, we give the definition of this
secondary multiplication in the next section.

\subsection{Secondary multiplication of the cohomology of a differential
graded algebra}

\label{sm} Let $C^{\ast}$ be a \emph{differential graded} algebra
(\emph{dg-algebra}) over a commutative hereditary ring $k$; that is, $C^{\ast
}$ is a graded $k$-algebra
\[
C=\oplus_{i\geq0}C^{i},\quad C^{i}C^{j}\subseteq C^{i+j},
\]
with a graded $k$-linear derivation (\emph{differential}) $d:C^{\ast
}\rightarrow C^{\ast}$ of degree 1
\[
B^{i+1}=d(C^{i})\subseteq C^{i+1},\quad d(xy)=d(x)y+(-1)^{\left|  x\right|
}xd(y),
\]
such that $d^{2}=0$. Here the formula is given for homogeneous $x,y$ and
$\left|  x\right|  $ is the degree of $x$.

Since $\operatorname*{Ker}(d)$ is a subalgebra in $C^{\ast}$ (subalgebra of
\emph{cocycles}) and $\operatorname*{Im}(d)$ is an ideal in
$\operatorname*{Ker}(d)$ (ideal of \emph{coboundaries}), the cohomology
$H^{\ast}(C)=\operatorname*{Ker}(d)/\operatorname*{Im}(d)$ of the dg-algebra
$C^{\ast}$ has a natural $k$-algebra structure. Here we consider dg-algebras
that arise as the dual of a free chain complex $(C_{\ast},\partial)$. The
following is easily checked (cf. \cite{Rotman} p.386). (We write
$B_{i-1}=\partial(C_{i}).)$

\begin{technical lemma}
\label{tech lemma}Suppose that $C_{\ast}$ is a free chain complex with
$H_{\ast}(C)$ a projective $k$-module. Then for each $n$ there is a direct sum
decomposition%
\[
C_{n}\cong H_{n}(C)\oplus B_{n}\oplus B_{n-1},
\]
which induces by duality a direct sum decomposition%
\begin{align*}
C^{n} &  \cong H_{n}(C)^{\vee}\oplus B_{n}^{\vee}\oplus B_{n-1}^{\vee}\\
&  \cong H^{n}(C)\oplus B^{n+1}\oplus B^{n}.\tag*{\qed}
\end{align*}\end{technical lemma}

It follows that we can choose graded $k$-linear sections $s:H^{\ast
}(C)\rightarrow\operatorname*{Ker}(d)$ and $q:\operatorname*{Im}(d)\rightarrow
C$ of the natural $k$-linear surjections $\pi:\operatorname*{Ker}%
(d)\rightarrow H^{\ast}(C)$ and $d:C\rightarrow\operatorname*{Im}(d)$
respectively. The situation may be summarized as
\[%
\begin{array}
[c]{ccccc}%
\operatorname*{Im}d &  &  &  & \\
\downarrow &  &  &  & \\
\operatorname*{Ker}d & \hookrightarrow &  C & \overset{\overset{q}%
{\curvearrowleft}}{\underset{d}{\twoheadrightarrow}} & \operatorname*{Im}%
d[-1].\\
{}^{\pi}\downarrow\uparrow^{s} &  &  &  & \\
H^{\ast}(C) &  &  &  &
\end{array}
\]

The map $\pi:\operatorname*{Ker}(d)\rightarrow H^{\ast}(C)$ is a $k$-algebra
homomorphism, hence the difference
\[
s(x)s(y)-s(xy)
\]
has zero image under $\pi$:
\[
\pi(s(x)s(y)-s(xy))=\pi(s(x))\pi(s(y))-\pi(s(xy))=xy-xy=0,
\]
that is, lies in $\operatorname*{Im}(d)$. Define $q(x,y)=q(s(x)s(y)-s(xy))\in
C$, so that
\[
dq(x,y)=s(x)s(y)-s(xy).
\]
(Observe that $\left|  q(x,y)\right|  =\left|  x\right|  +\left|  y\right|
-1.)$ Then the expression
\begin{equation}
\Theta(x,y,z)=(-1)^{\left|  x\right|  }s(x)q(y,z)-q(xy,z)+q(x,yz)-q(x,y)s(z)
\label{theta}%
\end{equation}
lies in $\operatorname*{Ker}(d)$. Indeed,
\begin{align*}
d\Theta(x,y,z)=  &  s(x)dq(y,z)-dq(xy,z)+dq(x,yz)-d(q(x,y))s(z)\\
=  &  s(x)(s(y)s(z)-s(yz))-(s(xy)s(z)-s(xyz))\\
&  +(s(x)s(yz)-s(xyz))-(s(x)s(y)-s(xy))s(z)\\
=  &  0.
\end{align*}
Define $\theta(x,y,z)=\pi(\Theta(x,y,z))\in H^{\ast}(C)$. Then
\[
\left|  \theta(x,y,z)\right|  =\left|  \Theta(x,y,z)\right|  =\left|
x\right|  +\left|  y\right|  +\left|  z\right|  -1,
\]
so $\theta:H^{\ast}(C)^{\otimes3}\rightarrow H^{\ast}(C)$ has degree $-1.$
This is a $k$-linear map $H^{\ast}(C)^{\otimes3}\rightarrow H^{\ast}(C)[1]$,
where $H^{\ast}(C)[1]$ is a \emph{shifted} graded module $H^{\ast}(C)$ such
that\break $(H^{\ast}(C)[1])^{m}=H^{m-1}(C)$. The twisted multiplication $x\ast
y=(-1)^{\left|  x\right|  }xy$ allows us to define a new $H^{\ast}%
(C)$-bimodule structure on $H^{\ast}(C)$ by setting the left module structure
to be twisted and the right module structure to be the ordinary one. We denote
this $H^{\ast}(C)$-bimodule by $\overline{H^{\ast}(C)}$.

\begin{proposition}
The map $\theta:H^{\ast}(C)^{\otimes3}\rightarrow H^{\ast}(C)$ defined above
is a Hoch\-schild 3-cocycle of degree -1 with respect to the graded algebra
$H^{\ast}(C)$.
\end{proposition}

\proof
We check the 3-cocycle property of $\theta:H^{\ast}(C)^{\otimes3}\rightarrow
H^{\ast}(C)$. Define the maps
\[
\delta_{s}:\operatorname*{Hom}(H^{\ast}(C)^{\otimes l},C)\rightarrow
\operatorname*{Hom}(H^{\ast}(C)^{\otimes l+1},C)
\]
as
\begin{align*}
&  (\delta_{s}\phi)(x_{1},...,x_{l+1})\\
=  &  (-1)^{\left|  x_{1}\right|  }s(x_{1})\phi(x_{2},...,x_{l+1})+\sum
_{i=1}^{l}(-1)^{i}\phi(x_{1},...,x_{i}x_{i+1},...,x_{l+1})\\
&  +(-1)^{l+1}\phi(x_{1},...,x_{l})s(x_{l+1}).
\end{align*}
The restrictions of these maps to $\operatorname*{Hom}(H^{\ast}(C)^{\otimes
l},\,\operatorname*{Ker}(d))$ are compatible with Hochschild differentials on
$\operatorname*{Hom}(H^{\ast}(C)^{\otimes l},\,\overline{H^{\ast}(C)}),$ so
that the diagram
\begin{equation}%
\begin{array}
[c]{ccccc}%
\overset{\delta_{s}}{\longrightarrow} & \operatorname*{Hom}(H^{\ast
}(C)^{\otimes l},\,\operatorname*{Ker}(d)) & \overset{\delta_{s}%
}{\longrightarrow} & \operatorname*{Hom}(H^{\ast}(C)^{\otimes l+1}%
,\,\operatorname*{Ker}(d)) & \overset{\delta_{s}}{\longrightarrow}\\
& \downarrow^{\pi\circ} &  & \downarrow^{\pi\circ} & \\
\overset{\delta}{\longrightarrow} & \operatorname*{Hom}(H^{\ast}(C)^{\otimes
l},\,\overline{H^{\ast}(C)}) & \overset{\delta}{\longrightarrow} &
\operatorname*{Hom}(H^{\ast}(C)^{\otimes l+1},\,\overline{H^{\ast}(C)}) &
\overset{\delta}{\longrightarrow}%
\end{array}
\label{Hochschild diagram}%
\end{equation}
is commutative. We shall show that
\[
\delta_{s}\Theta(x,y,z,w)=(-1)^{\left|  x\right|  +\left|  y\right|
}d(q(x,y)q(z,w)),
\]
which implies that $\theta=\pi\Theta$ is a 3-cocycle.

Since $s$ is not a ring homomorphism in general, the map $\delta_{s}$ does not
satisfy the condition $\delta_{s}\delta_{s}=0$. Instead of this, upon writing
$\delta_{s}=\sum_{i=0}^{l+1}(-1)^{i}\delta_{s}^{i}$ as usual, we note that
\[
\delta_{s}^{i}\delta_{s}^{j}=\delta_{s}^{j+1}\delta_{s}^{i}\qquad i\leq
j,\ (i,j)\neq(0,0),(l+1,l+1).
\]
So the Differential Lemma immediately implies that
\[
\delta_{s}\delta_{s}=(\delta_{s}^{0}\delta_{s}^{0}-\delta_{s}^{1}\delta
_{s}^{0})+(\delta_{s}^{l+1}\delta_{s}^{l+1}-\delta_{s}^{l+2}\delta_{s}%
^{l+1}),
\]
in other words,
\begin{align*}
&  \delta_{s}\delta_{s}\phi(x_{1},...,x_{l+2})\\
=  &  (-1)^{\left|  x_{1}\right|  +\left|  x_{2}\right|  }dq(x_{1},x_{2}%
)\phi(x_{3},...,x_{l+2})-\phi(x_{1},...,x_{l})dq(x_{l+1},x_{l+2}).
\end{align*}
Now, by definition, $\Theta=\delta_{s}q.$ Therefore
$$\delta_{s}\Theta(x,y,z,w)=(-1)^{\left|  x\right|  +\left|  y\right|
}d(q(x,y)q(z,w)).\eqno{\qed}
$$

\begin{proposition}
\label{choice of s q}The cohomology class $[\theta]\in HH^{3}(H^{\ast
}(C),\overline{H^{\ast}(C)}[1])$ is independent of the choices of $s$ and $q$.
\end{proposition}

\proof
For another choice $q^{\prime}$ of the section $q:\operatorname*{Im}%
(d)\rightarrow C,$ the difference $q^{\prime}-q$ is a map $\operatorname*{Im}%
(d)\rightarrow\operatorname*{Ker}(d)$. Hence we are in the situation of
diagram (\ref{Hochschild diagram}), and the cocycles $\theta^{\prime}$ and
$\theta$ differ by the coboundary $\delta(q^{\prime}-q)$.

For another choice $s^{\prime}$ of the section $s:H^{\ast}(C)\rightarrow C,$
the difference $s^{\prime}-s$ is a map $H^{\ast}(C)\rightarrow
\operatorname*{Im}(d)$. So we can write $s^{\prime}-s=da$ for some $a:H^{\ast
}(C)\rightarrow C$ of degree $-1$. Now because $ds=dd=0$, we have
\[
d((-1)^{\left|  x\right|  }s(x)a(y)+a(x)s(y)+a(x)da(y))=s^{\prime}%
(x)s^{\prime}(y)-s(x)s(y).
\]
Therefore
\begin{align*}
d(q^{\prime}(x,y)-q(x,y))  &  =s^{\prime}(x)s^{\prime}(y)-s^{\prime
}(xy)-(s(x)s(y)-s(xy))\\
&  =db(x,y),
\end{align*}
where
\[
b(x,y)=(-1)^{\left|  x\right|  }s(x)a(y)-a(xy)+a(x)s(y)+a(x)da(y).
\]
This means that the difference $q^{\prime}(x,y)-q(x,y)-b(x,y)$ lies in
$\operatorname*{Ker}(d)$ and, as we saw before, does not affect the cohomology
class. Hence we can suppose that $q^{\prime}(x,y)=q(x,y)+b(x,y)$. Then the
difference of cocycles
\begin{align*}
&  \Theta^{\prime}(x,y,z)-\Theta(x,y,z)\\
=\quad &  (-1)^{\left|  x\right|  }s^{\prime}(x)q^{\prime}(y,z)-q^{\prime
}(xy,z)+q^{\prime}(x,yz)-q^{\prime}(x,y)s^{\prime}(z)\\
-  &  ((-1)^{\left|  x\right|  }s(x)q(y,z)-q(xy,z)+q(x,yz)-q(x,y)s(z))
\end{align*}
can be rewritten as
\begin{align*}
&  (-1)^{\left|  x\right|  }da(x)(q(y,z)+b(y,z))-(q(x,y)+b(x,y))da(z)\\
+  &  (-1)^{\left|  x\right|  }s(x)b(y,z)-b(xy,z)+b(x,yz)-b(x,y)s(z).
\end{align*}
Using congruences modulo $\operatorname*{Im}(d),$
\begin{align*}
(-1)^{\left|  x\right|  }da(x)q(y,z)  &  \equiv a(x)dq(y,z)\\
&  =a(x)(s(y)s(z)-s(yz)),
\end{align*}%
\[
q(x,y)da(z)\equiv(-1)^{\left|  xy\right|  }(s(x)s(y)-s(xy))a(z),
\]
and the definition of $b,$ we can, after cancellation of similar terms,
rewrite the expression as
\[
(-1)^{\left|  x\right|  }d((-1)^{\left|  y\right|  }%
a(x)s(y)a(z)-a(x)a(yz)+a(x)a(y)s^{\prime}(z)),
\]
which lies in $\operatorname*{Im}(d),$ completing the proof.
\endproof

\remark
\label{mass}Connection with Massey triple product

\noindent This connection admits an abstract algebraic setting. For any triple
of elements $x,y,z\in R$ of an associative ring $R$ satisfying the conditions
$xy=yz=0$ (\emph{Massey triple}), we can define a homomorphism
\[
HH^{3}(R,M)\rightarrow M/(xM+Mz),
\]
by sending the class of a $3$-cocycle $\alpha$ to the image of its value
$\alpha(x,y,z)$ in $M/(xM+Mz)$. Indeed, the value of any coboundary
\begin{align}
\delta a(x,y,z)  &  =xa(y,z)-a(xy,z)+a(x,yz)-a(x,y)z\nonumber\\
&  =xa(y,z)-a(x,y)z
\end{align}
represents zero in $M/(xM+Mz)$.

Now, for any Massey triple $x,y,z\in H^{\ast}(C),$ the expression
$\Theta(x,y,z)$ takes the form
\[
(-1)^{\left|  x\right|  }s(x)q(y,z)-q(x,y)s(z),
\]
which coincides with the usual definition of $\left\langle x,y,z\right\rangle
$ (\emph{Massey triple product}) \cite{mas}. Note that the value of
$\left\langle x,y,z\right\rangle $ in the quotient $H^{\ast}(C)/(xH^{\ast
}(C)+H^{\ast}(C)z)$ is canonically defined, that is, does not depend on the
choice of $s$ and $q$. This property can be deduced from the canonical nature
of the definition of the corresponding Hochschild cohomology class
(Proposition \ref{choice of s q}).

Let $C^{\ast}(X;\,k)$ be the dg-algebra (under cup-product) of singular
cochains of a topological space $X$. It is the dual of the singular chain
complex $C_{\ast}(X;\,k)$, a free graded $k$-module. If the graded $k$-module
$H_{\ast}(X;\,k)$ is also projective, then (\ref{tech lemma}) applies, and we
have a linear section
\[
q:B^{\ast}(X;\,k)=\operatorname*{Im}(d)\rightarrow C^{\ast}(X;\,k)
\]
of the surjection $d:C^{\ast}(X;\,k)\rightarrow B^{\ast}(X;\,k),$ and section
$s:H^{\ast}(X;\,k)\rightarrow\mathrm{Ker}(d)$ of the surjection $\mathrm{Ker}%
(d)\rightarrow H^{\ast}(X;\,k).$ Hence the secondary multiplication class
\[
\left[  \theta(X,k)\right]  \in HH_{\operatorname*{gr}}^{3}(H^{\ast
}(X;\,k),\,\overline{H^{\ast}(X;\,k)}[1])
\]
is defined.

\section{Splitting of the multiplication map}

\label{smm}

For a morphism $f:K\rightarrow L$ of cochain complexes, we use the mapping
cone $\operatorname*{cone}(f),$ which fits into a short exact sequence of
complexes
\begin{equation}
0\rightarrow L\rightarrow\operatorname*{cone}(f)\rightarrow K[-1]\rightarrow0.
\label{et}%
\end{equation}
In particular, we have a long exact sequence of cohomology groups
\begin{equation}
\cdots\rightarrow H^{n}(L)\rightarrow H^{n}(\operatorname*{cone}%
(f))\rightarrow H^{n+1}(K)\overset{f}{\rightarrow}H^{n+1}(L)\rightarrow\cdots.
\label{les}%
\end{equation}
Recall that $\operatorname*{cone}(f)$ is defined by equipping the direct sum
$L\oplus K[-1]$ with the differential $D(x,y)=(d(x)+f(y),-d(y))$.

A (right) \emph{dg-module} over a dg-algebra $C$ is a complex $M$ with a
module structure over the algebra $C$
\begin{equation}
M\otimes C\rightarrow M,\quad\ m\otimes x\mapsto mx, \label{mst}%
\end{equation}
such that
\begin{equation}
d(mx)=d(m)x+(-1)^{\left|  m\right|  }md(x). \label{mlp}%
\end{equation}
In other words, a module over a dg-algebra is a module over an algebra such
that the module-structure map (\ref{mst}) is a homomorphism of complexes. A
\emph{morphism} or \emph{$C$-linear map} of dg-modules over the dg-algebra $C$
is a morphism of complexes that commutes with the module structure maps. The
following is proved by straightforward computation.

\begin{proposition}
Let $f:M\rightarrow N$ be a morphism of right dg-modules over a dg-algebra
$C$. Then the mapping cone $\operatorname*{cone}(f)$ has a natural dg-module
structure over $C$ such that all maps of the exact sequence (\ref{et}) are $C$-linear.\qed
\end{proposition}

The dg-module structure over the dg-algebra $C$ on the complex $M$ induces the
$H^{\ast}(C)$-module structure on its cohomology $H^{\ast}(M)$; any $C$-linear
map between dg-modules induces an $H^{\ast}(C)$-linear map between their
cohomology modules. In particular, for the $C$-linear map $f:M\rightarrow N$
between dg-modules, the sequence (\ref{les}) is an exact sequence of $H^{\ast
}(C)$-modules.

Now we examine some special cases of $C$-linear map. Note that the dg-algebra
$C$ can be considered as a module over itself. For any cocycle $z\in
\operatorname*{Ker}(d)\subset C,$ left multiplication by $z$ is a morphism of
complexes
\begin{equation}
l_{z}:C[\left|  z\right|  ]\rightarrow C, \label{mm}%
\end{equation}
which is also $C$-linear with respect to the natural right $C$-module
structure on $C$.

The map $H^{\ast}(C)[\left|  z\right|  ]\rightarrow H^{\ast}(C)$ induced by
the left multiplication map $l_{z}$ is left multiplication by the class
$c=\pi(z)\in H^{\ast}(C)$ of $z$. Thus in the case when$\ f=l_{s(c)}$ for some
$c\in H^{\ast}(C)$ we can rewrite the long exact sequence (\ref{les}) as an
extension of $H^{\ast}(C)$-modules
\begin{equation}
0\rightarrow H^{\ast}(C)/(cH^{\ast}(C)[\left|  c\right|  ])\rightarrow
H^{\ast}(\operatorname*{cone}(l_{s(c)}))\rightarrow\operatorname*{Ann}%
\nolimits_{H^{\ast}(C)}(c)[\left|  c\right|  -1]\rightarrow0, \label{ses}%
\end{equation}
where $\operatorname*{Ann}\nolimits_{H^{\ast}(C)}(c)=\{x\in H^{\ast}(C)\mid
cx=0\}$ is the right annihilator of $c$.

We shall describe the class of this extension in the group
\[
\mathcal{E}:=\operatorname*{Ext}\nolimits_{H^{\ast}(C)}^{1}%
(\operatorname*{Ann}\nolimits_{H^{\ast}(C)}(c)[\left|  c\right|
-1],\,H^{\ast}(C)/(cH^{\ast}(C)[\left|  c\right|  ])).
\]
We first develop an abstract setting for the connection between the groups
$HH_{\operatorname*{gr}}^{3}(H^{\ast}(C),\,\overline{H^{\ast}(C)}[1])$ and
$\mathcal{E}$. The following is easily checked.

\begin{lemma}
$\operatorname*{Ker}(\operatorname*{Ext}\nolimits_{R}^{1}(N,M)\rightarrow
\operatorname*{Ext}\nolimits_{k}^{1}(N,M))$ coincides with the group of
$k$-linear maps $\beta:N\otimes_{k}R\rightarrow M$ satisfying
\begin{equation}
\beta(n,xy)=\beta(nx,y)+\beta(n,x)y\label{ecl}%
\end{equation}
modulo its subgroup of maps of the form
\begin{equation}
\beta(n,x)=b(nx)-b(n)x,\ \label{tcl}%
\end{equation}
for some $k$-linear $b:N\rightarrow M.$\qed\end{lemma}

Let $R$ be a $k$-algebra and $M$ an $R$-bimodule. For any $c\in R$ we define a
homomorphism
\begin{equation}
HH^{3}(R,M)\rightarrow\operatorname*{Ext}\nolimits_{R}^{1}(\operatorname*{Ann}%
\nolimits_{R}(c),\,M/cM).
\end{equation}
Let $\alpha\in ZH^{3}(R,M)$ be a cocycle of the standard Hochschild complex.
Define the map $\beta:\operatorname*{Ann}\nolimits_{R}(c)\otimes R\rightarrow
M/cM$ by setting $\beta(x,y)$ to be the class of $\alpha(c,x,y)$ in $M/cM$.
The cocycle property
\[
c\alpha(x,y,z)-\alpha(cx,y,z)+\alpha(c,xy,z)-\alpha(c,x,yz)+\alpha(c,x,y)z=0
\]
implies the equation (\ref{ecl})
\[
\beta(x,yz)=\beta(xy,z)+\beta(x,y)z,
\]
for any $x\in\operatorname*{Ann}_{R}(c)$ and $y,z\in R$. Thus, by the lemma
above, $\beta$ defines the structure of an $R$-extension on
$\operatorname*{Ann}{}_{R}(c)\oplus M/cM$. For a coboundary $\alpha=d(a)$
\[
\alpha(c,x,y)=ca(x,y)-a(cx,y)+a(c,xy)-a(c,x)y
\]
the corresponding $\beta$ takes the form $\beta(x,y)=b(xy)-b(x)y$ where
$b(x)=a(c,x)$.

Now let $R$ be a graded $k$-algebra, $M$ a graded $R$-bimodule and $c\in R$
homogeneous of degree $\left|  c\right|  $. For $\alpha$ of degree zero, the
degree of its corresponding $\beta$ is $\left|  c\right|  $. Hence we have a
map
\[
HH_{\operatorname*{gr}}^{3}(R,M)\rightarrow\operatorname*{Ext}\nolimits_{R}%
^{1}(\operatorname*{Ann}{}_{R}(c)[\left|  c\right|  ],\,M/cM).
\]
In particular, when $R=H^{\ast}(C)$ and $M=\overline{H^{\ast}(C)}[1],$ we have
the map
\begin{equation}
HH_{\operatorname*{gr}}^{3}(H^{\ast}(C),\,\overline{H^{\ast}(C)}%
[1])\rightarrow\mathcal{E}. \label{hme}%
\end{equation}

\begin{theorem}
\label{th} Let $C$, a dg-algebra over a hereditary ring, be the dual of a
projective chain complex with projective homology$.$ Let $c\in H^{\ast}(C),$
and let
\[
\eta\in\operatorname*{Ext}{}_{H^{\ast}(C)}^{1}(\operatorname*{Ann}{}_{H^{\ast
}(C)}(c)[\left|  c\right|  -1],\,H^{\ast}(C)/(cH^{\ast}(C)[\left|  c\right|
]))
\]
be the class of the extension (\ref{ses}). Then the map (\ref{hme}) sends the
Hochschild cohomology class $[\theta]$ to the class $\eta$.
\end{theorem}%

\proof
Note that the mapping cone $\operatorname*{cone}(l_{s(c)})$ of the
multiplication map (\ref{mm}) coincides with $C\oplus C[\left|  c\right|  -1]$
equiped with the differential
\[
D(x,y)=(d(x)+s(c)y,(-1)^{\left|  c\right|  -1}d(y)).
\]
We construct a section of the surjection
\begin{equation}
H^{\ast}(\operatorname*{cone}(l_{s(c)}))\rightarrow\operatorname*{Ann}%
{}_{H^{\ast}(C)}(c)[\left|  c\right|  -1] \label{sj}%
\end{equation}
from the short exact sequence (\ref{ses}), as follows. For $x\in
\operatorname*{Ann}{}_{H^{\ast}(C)}(c)$ the product $s(c)s(x)$ lies in
$B^{\ast}(C),$ so we can write $s(c)s(x)=dq(c,x)$ for a section $q$ of $d$ as
in (\ref{sm}). Then $(-q(c,x),s(x))$ is a cocycle of $\operatorname*{cone}%
(l_{s(c)}),$ because
\[
D(-q(c,x),s(x))=(-dq(c,x)+s(c)s(x),(-1)^{\left|  c\right|  -1}ds(x))=0.
\]
Obviously, the image of $(-q(c,x),s(x))$ under the natural map
$\operatorname*{cone}(l_{s(c)})\rightarrow C[\left|  c\right|  -1]$ is $s(x)$.
Thus the map $\sigma:\operatorname*{Ann}{}_{H^{\ast}(C)}(c)[\left|  c\right|
-1]\rightarrow H^{\ast}(\operatorname*{cone}(l_{s(c)}))$ sending $x$ to the
class of $(-q(c,x),s(x))$ in $H^{\ast}(\operatorname*{cone}(l_{s}(c)))$ is a
section of the surjection (\ref{sj}).

Now we calculate the expression $\sigma(x)y-\sigma(xy)$ which represents the
class $\eta$ of the extension (\ref{ses}) in $\mathcal{E}.$
\begin{align*}
\sigma(x)y-\sigma(xy)  &  =\pi((-q(c,x),\,s(x))s(y)-(-q(c,xy),s(xy)))\\
&  =\pi(-q(c,x)s(y)+q(c,xy),\,s(x)s(y)-s(xy))\\
&  =\pi(-q(c,x)s(y)+q(c,xy),\,dq(x,y)).
\end{align*}
Since $(0,dq(x,y))=((-1)^{\left|  c\right|  }s(c)q(x,y),0)+(-1)^{\left|
c\right|  -1}D(0,q(x,y)),$ we have
\[
\sigma(x)y-\sigma(xy)=\pi((-1)^{\left|  c\right|  }%
s(c)q(x,y)+q(c,xy)-q(c,x)s(y),0),
\]
which from (\ref{theta}) coincides with the image of $\theta(c,x,y)$ under the
map $H^{\ast}(C)\rightarrow H^{\ast}(\operatorname*{cone}(l_{s(C)}))$.
\endproof

\section{Cohomology of an orientable sphere bundle}

In order to apply the previous algebraic result (\ref{th}) to geometrical
situations, we use the following.

\begin{lemma}
Let $f:C_{1}\rightarrow C_{2}$ be a homomorphism of dg-algebras and $z\in
C_{1}$ be a cocycle whose class $c$ lies in the kernel $c\in
\operatorname*{Ker}H(f)$. Then the homomorphism $f$ of dg-algebras extends to
a (right) $C_{1}$-linear map of complexes $\operatorname*{cone}(l_{z}%
)\rightarrow C_{2}$.
\end{lemma}

\proof
Since the class $c$ of the cocycle $z\in C_{1}$ lies in the kernel of the
homomorphism $H^{\ast}(C_{1})\rightarrow H^{\ast}(C_{2}),$ we can find some
cochain $b\in C_{2}^{\left|  z\right|  -1}$ satisfying $db=f(z)$. Thus we can
define a map $\operatorname*{cone}(l_{z})\rightarrow C_{2}$ sending $(x,y)$ to
$f(x)+bf(y)$. This map is compatible with the differentials, for
\begin{align*}
d(f(x)+bf(y))  &  =df(x)+(db)f(y)+(-1)^{\left|  b\right|  }b(df(y))\\
&  =f(dx)+f(zy)+(-1)^{\left|  z\right|  -1}bf(dy)
\end{align*}
coincides with the image of $D(x,y)=(d(x)+zy,\,(-1)^{\left|  z\right|
-1}d(y)) $. Moreover, the map is obviously $C_{1}$-linear with respect to the
structure of a right dg-module over $C_{1}$ on $C_{2}$ given by the dg-algebra
homomorphism $f:C_{1}\rightarrow C_{2}$.\endproof

Now let $S^{m}\rightarrow X\overset{p_{0}}{\longrightarrow}B$ be an orientable
sphere bundle, with associated disc bundle $D^{m+1}\rightarrow E\overset
{p}{\longrightarrow}B$. Denote by $u\in H^{m+1}(E,X)$ its orientation class
and by $c\in H^{m+1}(B)$ its characteristic class. Choose some representative
$z\in Z^{m+1}(B)$ for $c$. Write $l_{z}:C^{\ast}(B)[m+1]\rightarrow C^{\ast
}(B)$ for the left multiplication map. The class $c$ lies in the kernel of the
homomorphism $p_{0}^{\ast}:H^{m+1}(B)\rightarrow H^{m+1}(X)$. Thus we can
apply the lemma and extend the homomorphism of dg-algebras $p_{0}^{\ast
}:C^{\ast}(B)\rightarrow C^{\ast}(X)$ to a $C^{\ast}(B)$-linear map of
complexes $\bar{p}:\operatorname*{cone}(l_{z})\rightarrow C^{\ast}(X)$.

\begin{proposition}
The map $\operatorname*{cone}(l_{z})\rightarrow C^{\ast}(X)$ induces an
$H^{\ast}(B)$-linear isomorphism in cohomologies, and the long exact sequence
(\ref{les}) is isomorphic under this map to the Gysin sequence of the sphere
bundle
\[%
\begin{tabular}
[c]{ccccccccc}%
$\longrightarrow$\hspace{-7pt}& $H^{\ast}(B)$ & $\longrightarrow$ & $H^{\ast
}(\operatorname*{cone}(l_{z}))$ & $\longrightarrow$ & $H^{\ast-m}(B)$ &
$\overset{l_{c}}{\longrightarrow}$ & $H^{\ast+1}(B)$ & $\hspace{-7pt}\longrightarrow$\\
& $\Vert$ &  & $\downarrow^{H(\bar{p})}$ &  & $\Vert$ &  & $\Vert$ & \\
$\longrightarrow$\hspace{-7pt} & $H^{\ast}(B)$ & $\longrightarrow$ & $H^{\ast}(X)$ &
$\longrightarrow$ & $H^{\ast-m}(B)$ & $\overset{l_{c}}{\longrightarrow}$ &
$H^{\ast+1}(B)$ &\hspace{-7pt} $\longrightarrow$%
\end{tabular}
\]
\end{proposition}

\proof
It is well-known that the Gysin sequence of a sphere bundle can be derived
from the long exact sequence of the pair $(E,X),$ using the natural
identification $H^{\ast}(B)\overset{p^{\ast}}{\longrightarrow}H^{\ast}(E)$ and
the Thom isomorphism $H^{\ast}(B)[m+1]\overset{l_{u}}{\longrightarrow}H^{\ast
}(E,X)$ given by multiplication by the orientation class $u\in H^{m+1}(E,X).$

We can establish this connection on the level of singular cochain complexes.
The short exact sequence $0\rightarrow C^{\ast}(E,X)\rightarrow C^{\ast
}(E)\rightarrow C^{\ast}(X)\rightarrow0,$ which gives the long exact sequence
of the pair $(E,X),$ can be extended to the diagram
\[%
\begin{tabular}
[c]{lllllllll}%
&  & $C^{\ast}(B)[m+1]$ & $\overset{l_{z}}{\longrightarrow}$ & $C^{\ast}(B)$ &
$\overset{p_{0}^{\ast}}{\longrightarrow}$ & $C^{\ast}(X)$ &  & \\
&  & $\downarrow^{l_{U}\circ p^{\ast}}$ &  & $\downarrow^{p^{\ast}}$ &  &
$\Vert$ &  & \\
$0$ & $\longrightarrow$ & $C^{\ast}(E,X)$ & $\longrightarrow$ & $C^{\ast}(E)$%
& $\overset{}{\longrightarrow}$ & $C^{\ast}(X)$ & $\longrightarrow$ & $0$%
\end{tabular}
\]
where the first vertical arrow is a composition of $p^{\ast}:C^{\ast
}(B)\rightarrow C^{\ast}(E)\ $with multiplication by a representative cocycle
$U\in C^{m+1}(E,X)$ of the orientation class $u\in H^{m+1}(E,X).$ The
right-hand square of the diagram is obviously commutative. We seek to choose
cocycles $z$ and $U$ such that the left-hand square of the diagram is
commutative up to homotopy. To achieve this, we exploit the following fact.
Denote by $s:B\rightarrow E$ the zero section of the associated disc bundle.
The composition $ps$\ is the identity, so $s^{\ast}p^{\ast}=\operatorname*{id}%
\nolimits_{C^{\ast}(B)}.$ Further, $sp$\ is homotopic to the identity, so we
have $p^{\ast}s^{\ast}=\operatorname*{id}_{C^{\ast}(E)}+dh+hd$ for some
$h:C^{\ast}(E)\rightarrow C^{\ast}(B)[1].$ Now we can start with a given $U$
and define $z$ to be $s^{\ast}U.$\ Then the homotopy that makes the left
square of the diagram commutative coincides, up to sign, with left
multiplication by $h(U)$ composed with $p^{\ast}$. Indeed,
\begin{align*}
p^{\ast}l_{z}  &  =l_{p^{\ast}z}\circ p^{\ast}\\
&  =l_{p^{\ast}s^{\ast}U}\circ p^{\ast}\\
&  =(l_{U}+l_{dh(U)})p^{\ast}\\
&  =l_{U}p^{\ast}+dl_{h(U)}p^{\ast}+(-1)^{m}l_{h(U)}p^{\ast}d.
\end{align*}
The upper row of the diagram is not exact. But we can extend it to the
morphism (up to homotopy) of two short exact sequences
\[%
\begin{tabular}
[c]{lllllllll}%
$0$ & $\longrightarrow$ & $C^{\ast}(B)[m+1]$ & $\longrightarrow$ &
$\operatorname*{cyl}(l_{z})$ & $\longrightarrow$ & $\operatorname*{cone}%
(l_{z})$ & $\longrightarrow$ & $0$\\
&  & $\Vert$ &  & $\downarrow^{\rho}$ &  & $\downarrow^{\bar{p}}$ &  & \\
&  & $C^{\ast}(B)[m+1]$ & $\overset{l_{z}}{\longrightarrow}$ & $C^{\ast}(B)$ &
$\overset{q^{\ast}}{\longrightarrow}$ & $C^{\ast}(X)$ &  & \\
&  & $\downarrow^{l_{U}\circ p^{\ast}}$ &  & $\downarrow^{p^{\ast}}$ &  &
$\Vert$ &  & \\
$0$ & $\longrightarrow$ & $C^{\ast}(E,X)$ & $\longrightarrow$ & $C^{\ast}(E)$%
& $\longrightarrow$ & $C^{\ast}(X)$ & $\longrightarrow$ & $0$%
\end{tabular}
\]
where the upper row is the short exact sequence of the cylinder
$\operatorname*{cyl}(l_{z})$ of the left multiplication $l_{z}$. The
quasi-isomorphism $\rho$ is left inverse to the inclusion of $C^{\ast}(B)$ in
$\operatorname*{cyl}(l_{z}).$ We conclude with the stock remark that, up to a
shift, the long exact sequence corresponding to the short exact sequence of
the cylinder $\operatorname*{cyl}(l_{z})$ coincides with the long exact
sequence corresponding to the short exact sequence of the cone.\endproof

Combining this with Theorem \ref{th}, we obtain the following.

\begin{theorem}
\label{mco} Let $S^{m}\rightarrow X\rightarrow B$ be an orientable sphere
bundle, where the homology $H_{\ast}(B;\,k)$ of the base is a projective
$k$-module. Let $c\in H^{m+1}(B)$ be its characteristic class and
\[
\eta\in\operatorname*{Ext}\nolimits_{H^{\ast}(B)}^{1}(\operatorname*{Ann}%
{}_{H^{\ast}(B)}(c)[m],\,H^{\ast}(B)/(cH^{\ast}(B)[m+1]))
\]
the class of the Gysin extension
\[
0\rightarrow H^{\ast}(B)/(cH^{\ast}(B)[m+1])\rightarrow H^{\ast}%
(X)\rightarrow\operatorname*{Ann}{}_{H^{\ast}(B)}(c)[m]\rightarrow0.
\]
Then the map (\ref{hme}) sends the Hochschild cohomology class of $\theta(B)$
to the class $\eta$.
\end{theorem}

\section{Applications and extensions}

Our information on the vanishing of the secondary multiplication class now
leads to two types of results.

\subsection{Rational coefficients and formal base}

In the rational case, we may apply to formal spaces, as follows.

A space $X$ is called (rationally) \emph{formal} if there is an $A_{\infty}%
$-algebra map \cite{kad1} $H^{\ast}(X,\mathbb{Q})\rightarrow C^{\ast
}(X,\mathbb{Q})$ inducing an isomorphism in cohomology. According to
\cite{kad1}, this is equivalent to triviality of the $A_{\infty}$-structure on
$H^{\ast}(X,\mathbb{Q})$. In particular (as can be seen directly from the
definition), the secondary multiplication class [$\theta(X,\mathbb{Q})]$ of a
formal space vanishes.

Rational homotopy theory has provided many examples of formal
spaces.\break Among them are simply-connected compact K\"{a}hler
manifolds \cite{dgms}, complex smooth algebraic varieties \cite{morg},
$(k-1)$-connected compact manifolds of dimension less than or equal to
$4k-2$ $($where $k\geq2)$ \cite{mill}, some flag varieties \cite{kum},
and some loop spaces \cite{dup,par}.

\begin{theorem}
Let $S^{m}\rightarrow X\rightarrow B$ be an orientable sphere bundle having
base $B$ a formal space, and characteristic class $c\in H^{m+1}(B)$. Then the
Gysin sequence with rational coefficients
\[
0\rightarrow H^{\ast}(B)/(cH^{\ast}(B)[m+1])\rightarrow H^{\ast}%
(X)\rightarrow\operatorname*{Ann}\nolimits_{H^{\ast}(B)}(c)[m]\rightarrow0
\]
splits as a sequence of right $H^{\ast}(B)$-modules.
\end{theorem}

\subsection{Integer coefficients and torus base\label{torus}}

Of course, rational formality does not in general imply formality over integer
or finite coefficients (see, for example \cite{eke}). However, when the base
is a torus $T$ and one takes integer coefficients, because the homology is
free abelian one can take advantage of the technical lemma (\ref{tech lemma}).

Here we briefly indicate our calculation of the secondary multiplication class
$\left[  \theta(T,\mathbb{Z})\right]  $ of an $n$-torus $T=\mathbf{T}%
^{n}=(A\otimes\mathbb{R})/(A\otimes1),$ using the quasi-isomorphism between
the dg-algebra of singular cochains $C^{\ast}(T)$ and the standard cochain
dg-algebra (bar construction) $C^{\ast}(A)$ of the free abelian group $A$ of
rank $n$.

Since the graded $\mathbb{Z}$-module $C_{\ast}(A)$ is free and also $H_{\ast
}(A)$ is free abelian, by (\ref{tech lemma}) we have a linear section
$q:B^{\ast}(A)=\operatorname*{Im}(d)\rightarrow C^{\ast}(A)$ of the surjection
$d:C^{\ast}(A)\rightarrow B^{\ast}(A)$ and a linear section $s:H^{\ast
}(A;\,k)\rightarrow\operatorname*{Ker}(d)$ of the natural projection. So, by
Proposition \ref{choice of s q}, the class%
\[
\left[  \theta(A)\right]  \in HH_{\operatorname*{gr}}^{3}(H^{\ast
}(A),\,\overline{H^{\ast}(A)}[1])
\]
is defined, and does not depend on the choices of sections $q$ and $s$ made above.

The cohomology $H^{\ast}(A)$ of the free abelian group $A$ coincides with the
algebra $\operatorname*{Hom}(\Lambda^{\ast}A,\,\mathbb{Z})$ of alternating
polylinear maps
\[
\operatorname*{Hom}(\Lambda^{l}A,\,\mathbb{Z})=\left\{  x:A^{l}\rightarrow
\mathbb{Z}\mid\ x(a_{1},...,a_{l})=0\
\mbox{whenever}%
\ a_{i}=a_{j}\ i\neq j\right\}  .
\]
We can consider $\operatorname*{Hom}(\Lambda^{\ast}A,\,\mathbb{Z})$ as the
exterior algebra $\Lambda^{\ast}V$ of the $\mathbb{Z}$-module $V=A^{\vee
}=\operatorname*{Hom}(A,\,\mathbb{Z}).$

To describe the graded Hochschild cohomology $HH_{\operatorname*{gr}}%
^{l}(R,\bar{R}[1])$ of the exterior algebra $R=\Lambda^{\ast}V$, we start with
construction of a map
\begin{equation}
HH_{\operatorname*{gr}}^{l}(\Lambda^{\ast}V,\,\overline{\Lambda^{\ast}%
V}[1])\rightarrow\operatorname*{Hom}(S^{l}V,\,\Lambda^{l-1}V) \label{dc}%
\end{equation}
for each $l,$ where $S^{l}V=(V^{\otimes l})_{S_{l}}$ is the module of
coinvariants of the natural action of the symmetric group $S_{l}$ on
$V^{\otimes l}$. So, for an abelian group $M,$ $\operatorname*{Hom}%
(S^{l}V,\,M)$ can be identified with the abelian group of symmetric polylinear
maps from $V$ to $M$. Now the above map is obtained from a homomorphism of
complexes induced from symmetrization, and can thereby be reduced to a
homomorphism involving a Koszul-type complex as in \cite{CartanEilenberg}%
(IX.6). Then, using induction on the rank of $A,$ it can be shown to be a
monomorphism. Explicit computation of the image of $\left[  \theta(A)\right]
$ in $\operatorname*{Hom}(S^{3}V,\,\Lambda^{2}V)$, exploiting symmetrization,
gives this image zero. Thus we obtain the desired conclusion.

\begin{proposition}
The class $\left[  \theta(A)\right]  \in HH_{\operatorname*{gr}}^{3}(H^{\ast
}(A),\,\overline{H^{\ast}(A)}[1])$ is trivial.\qed\end{proposition}

As a direct consequence of the quasi-isomorphism of dg-algebras between
$C^{\ast}(T)$ and $C^{\ast}(A),$ we obtain the vanishing of the topological obstruction.

\begin{corollary}
\label{clt} For a torus $T,$ the class of $\theta(T)$ in
$HH_{\operatorname*{gr}}^{3}(H^{\ast}(T),\,\overline{H^{\ast}(T)}[1])$ is
trivial.\qed\end{corollary}

The vanishing of the class of $\theta(T)$ gives the desired splitting.

\begin{theorem}
\label{tb} Let $S^{m}\rightarrow X\rightarrow T$ be an orientable sphere
bundle with torus base $T$ and characteristic class $c\in H^{m+1}(T)$. Then
the Gysin sequence with integer coefficients%
\[
0\rightarrow H^{\ast}(T)/(cH^{\ast}(T)[m+1])\rightarrow H^{\ast}%
(X)\rightarrow\operatorname*{Ann}\nolimits_{H^{\ast}(T)}(c)[m]\rightarrow0
\]
splits as a sequence of right $H^{\ast}(T)$-modules.\qed\end{theorem}

\subsection{Extensions}

\paragraph{Fibrations}

Although we have, for simplicity, presented our results in terms of sphere
bundles, they are also applicable to the generalized Gysin sequence of an
orientable fibration with fibre $F$ a cohomology $m$-sphere \cite{Spanier}%
(9.5.2). In the more general setting, $E$ may be obtained by the fibrewise
cone construction on the fibration $F\rightarrow X\overset{p_{0}%
}{\longrightarrow}B$ \cite{DrorFarjoun}(1.F). In order to repeat our argument
on the singular cochain complexes for the pair $(E,X),$ we need to show that
the Gysin sequence is again the cohomology exact sequence of the pair $(E,X).$
Chasing these two sequences, one observes that, because $E\rightarrow B$ is a
homotopy equivalence and by assumption $H^{m+1}(B)$ is $k$-torsion-free, both
of the groups $H^{0}(B)$ and $H^{m+1}(E,X)$ are isomorphic to the kernel of
$p_{0}^{\ast}:H^{m+1}(B)\rightarrow H^{m+1}(X)$ when the characteristic class
in $H^{m+1}(B)$ is nonzero, and to the cokernel of $p_{0}^{\ast}%
:H^{m}(B)\rightarrow H^{m}(X)$ otherwise. Corresponding to a generator of
$H^{0}(B)$ we therefore obtain an orientation class $u$ in $H^{m+1}(E,X).$
Since the two exact sequences coincide when $B$ is restricted to any point,
there results a cohomology extension of the fibre as in \cite{Spanier}(5.7),
and so a generalized Thom isomorphism $H^{i}(B)\rightarrow H^{i+m+1}(E,X)$ as required.

\paragraph{Bimodules}

Skew-commutativity of singular cohomology implies that our\break main results
(Theorems \ref{mco}, \ref{tb}) are valid in the context of bimodules, and not
only for right modules. Corresponding results of Section \ref{smm} can also be
proven for this context; however, the proofs are much more involved and need
some additional conditions for considering dg-algebras, such as stronger
homotopy commutativity. That is the reason we restricted ourselves there to
the case of right modules.

\Addresses\recd
 
\end{document}

%% file: agtout.tex
\def\ifplaintex{\expandafter\ifx\csname documentclass\endcsname\relax}

\def\gtp{{\mathsurround=0pt\it $\cal G\mskip-2mu$eometry \&\ 
$\cal T\!\!$opology $\cal P\!$ublications}}  

\def\recd{{\small Received:\qua\receiveddate\ifx\reviseddate\relax
\else\qquad Revised:\qua\reviseddate\fi\par}} 

\def\lognumber#1{\def\thelognumber{#1}}
\def\volumenumber#1{\def\thevolumenumber{#1}}
\def\volumeyear#1{\def\thevolumeyear{#1}}
\def\papernumber#1{\def\thepapernumber{#1}}
\def\pagenumbers#1#2{\def\startpage{#1}\def\finishpage{#2}}
\def\published#1{\def\publishdate{#1}}

\def\received#1{\def\receiveddate{#1}}
\def\revised#1{\def\reviseddate{#1}}
\def\accepted#1{\def\accepteddate{#1}}

\def\asciiaddress#1{\def\theasciiaddress{#1}}

\long\def\asciiabstract#1{\long\def\theasciiabstract{#1}}
\def\asciikeywords#1{\def\theasciikeywords{#1}}

\let\\\par\let\thelognumber\relax\let\thevolumenumber\relax
\let\thepapernumber\relax\let\thevolumeyear\relax\let\startpage\relax
\let\finishpage\relax\let\publishdate\relax\let\receiveddate\relax
\let\reviseddate\relax\let\accepteddate\relax\let\theasciititle\relax
\let\theasciiauthors\relax\let\theasciiaddress\relax
\let\theasciiabstract\relax\let\theasciikeywords\relax

\let\theasciiemail\relax

\ifplaintex
\font\logobig=cmssbx10 scaled 3836
\font\logomed=cmssbx10 scaled 2557
\else
\font\logobig=cmssbx10 scaled 4200
\font\logomed=cmssbx10 scaled 2800
\fi

\long\def\makeagttitle{   
\count0=\startpage
\agt\hfill      
\hbox to 45truept{\vbox to 0pt{\vglue -13truept{\logomed A\kern -.37em{\logobig 
T}\kern -.38em G}\vss}\hss}
\break
{\small Volume \thevolumenumber\ (\thevolumeyear)
\startpage--\finishpage\nl
Published: \publishdate}

\vglue .25truein

{\parskip=0pt\leftskip 0pt plus
1fil\def\\{\par\smallskip}{\Large\bf\thetitle}\par\medskip} \vglue
0.05truein

{\parskip=0pt\leftskip 0pt plus 1fil\def\\{\par}{\sc\theauthors}
\par\medskip}%
 
\vglue 0.03truein 

{\small\leftskip 25truept\rightskip 25truept{\bf Abstract}\stdspace\theabstract

{\bf AMS Classification}\stdspace\theprimaryclass
\ifx\thesecondaryclass\relax\else; \thesecondaryclass\fi\par
{\bf Keywords}\stdspace \thekeywords\par}\vglue 7truept

}   

\ifplaintex
\font\phead=cmsl9 scaled 950
\font\pnum=cmbx10 scaled 913
\font\pfoot=cmsl9 scaled 950
\headline{\vbox to 0pt{\vskip -4.5mm\line{\small\phead\ifnum
\count0=\startpage ISSN 1472-2739 (on-line) 1472-2747 (printed)
\hfill {\pnum\folio}\else\ifodd\count0\def\\{ }%
\ifx\theshorttitle\relax\thetitle\else\theshorttitle\fi\hfill{\pnum\folio}
\else\def\\{ and }{\pnum\folio}\hfill\ifx\theshortauthors\relax\theauthors
\else\theshortauthors\fi\fi\fi}\vss}}
\footline{\vbox to 0pt{\vglue 0mm\line{\small\pfoot\ifnum\count0=\startpage
\copyright\ \gtp\hfill\else
\agt, Volume \thevolumenumber\ (\thevolumeyear)\hfill\fi}\vss}}
\else
\font=cmsl9 scaled 1050
\font\lnum=cmbx10 
\font=cmsl9 scaled 1050
\makeatletter
\def\@oddhead{{\small\lhead\ifnum\count0=\startpage ISSN 1472-2739 
(on-line) 1472-2747 (printed)\hfill {\lnum\number\count0}\else\ifodd\count0
\def\\{ }\ifx\theshorttitle\relax \thetitle \else\theshorttitle\fi\hfill
{\lnum\number\count0}\else\def\\{ and }{\lnum\number\count0}
\hfill\ifx\theshortauthors\relax 
\theauthors\else\theshortauthors\fi\fi\fi}}\def\@evenhead{\@oddhead}
\def\@oddfoot{\small\lfoot\ifnum\count0=\startpage\copyright\ \gtp\hfill\else
\agt, Volume \thevolumenumber\ (\thevolumeyear)\hfill\fi}
\def\@evenfoot{\@oddfoot}
\makeatother
\fi
\let\maketitlepage\makeagttitle
\let\makeshorttitle\maketitlepage
\let\maketitle\maketitlepage

\newwrite\gtoutfile
\long\gdef\makeheadfile{  
{\def\\{, }\def\s{ }
\immediate\openout\gtoutfile head.xxx
\immediate\write\gtoutfile{To: math@arxiv.org}
\immediate\write\gtoutfile{Subject: put OR rep NNNNN:ppppp}
\immediate\write\gtoutfile{--text follows this line--}
\immediate\write\gtoutfile{Proxy-for: \ifx\theasciiauthors\relax
\theauthors\else\theasciiauthors\fi\s<\ifx\theasciiemail\relax\theemail\else\theasciiemail\fi>}
\immediate\write\gtoutfile{\noexpand\\}
\immediate\write\gtoutfile{Authors: \ifx\theasciiauthors\relax
\theauthors\else\theasciiauthors\fi}
{\def\\{ }\immediate\write\gtoutfile{Title: \ifx\theasciititle\relax
\thetitle\else\theasciititle\fi}}
\immediate\write\gtoutfile{Subj-class: GT or SG, GR etc}
\immediate\write\gtoutfile{MSC-class: \theprimaryclass\ifx\thesecondaryclass\relax\else, \thesecondaryclass\fi}
\immediate\write\gtoutfile{Journal-ref: Algebr. Geom. Topol. \thevolumenumber\s
(\thevolumeyear) \startpage-\finishpage}
\immediate\write\gtoutfile{Comments: Published by Algebraic and
Geometric Topology at}
\immediate\write\gtoutfile{\s\s\s  http://www.maths.warwick.ac.uk/agt/AGTVol\thevolumenumber/agt-\thevolumenumber-\thepapernumber.abs.html}
\immediate\write\gtoutfile{\noexpand\\}
\immediate\write\gtoutfile{}
\ifx\theasciiabstract\relax
\immediate\write\gtoutfile{\theabstract}\else
\immediate\write\gtoutfile{\theasciiabstract}\fi
\immediate\write\gtoutfile{}
\immediate\write\gtoutfile{\noexpand\\}
\immediate\write\gtoutfile{}
\immediate\closeout\gtoutfile}}  

\def\maketitlepage{\makeagttitle\makeheadfile}
\let\makeshorttitle\maketitlepage
\let\maketitle\maketitlepage

\def\ifplaintex{\expandafter\ifx\csname documentclass\endcsname\relax}

\def\gtp{{\mathsurround=0pt\it $\cal G\mskip-2mu$eometry \&\ 
$\cal T\!\!$opology $\cal P\!$ublications}}  

\def\recd{{\small Received:\qua\receiveddate\ifx\reviseddate\relax
\else\qquad Revised:\qua\reviseddate\fi\par}} 

\def\lognumber#1{\def\thelognumber{#1}}
\def\volumenumber#1{\def\thevolumenumber{#1}}
\def\volumeyear#1{\def\thevolumeyear{#1}}
\def\papernumber#1{\def\thepapernumber{#1}}
\def\pagenumbers#1#2{\def\startpage{#1}\def\finishpage{#2}}
\def\published#1{\def\publishdate{#1}}

\def\received#1{\def\receiveddate{#1}}
\def\revised#1{\def\reviseddate{#1}}
\def\accepted#1{\def\accepteddate{#1}}

\def\asciiaddress#1{\def\theasciiaddress{#1}}

\long\def\asciiabstract#1{\long\def\theasciiabstract{#1}}
\def\asciikeywords#1{\def\theasciikeywords{#1}}

\let\\\par\let\thelognumber\relax\let\thevolumenumber\relax
\let\thepapernumber\relax\let\thevolumeyear\relax\let\startpage\relax
\let\finishpage\relax\let\publishdate\relax\let\receiveddate\relax
\let\reviseddate\relax\let\accepteddate\relax\let\theasciititle\relax
\let\theasciiauthors\relax\let\theasciiaddress\relax
\let\theasciiabstract\relax\let\theasciikeywords\relax

\let\theasciiemail\relax

\ifplaintex
\font\logobig=cmssbx10 scaled 3836
\font\logomed=cmssbx10 scaled 2557
\else
\font\logobig=cmssbx10 scaled 4200
\font\logomed=cmssbx10 scaled 2800
\fi

\long\def\makeagttitle{   
\count0=\startpage
\agt\hfill      
\hbox to 45truept{\vbox to 0pt{\vglue -13truept{\logomed A\kern -.37em{\logobig 
T}\kern -.38em G}\vss}\hss}
\break
{\small Volume \thevolumenumber\ (\thevolumeyear)
\startpage--\finishpage\nl
Published: \publishdate}

\vglue .25truein

{\parskip=0pt\leftskip 0pt plus
1fil\def\\{\par\smallskip}{\Large\bf\thetitle}\par\medskip} \vglue
0.05truein

{\parskip=0pt\leftskip 0pt plus 1fil\def\\{\par}{\sc\theauthors}
\par\medskip}%
 
\vglue 0.03truein 

{\small\leftskip 25truept\rightskip 25truept{\bf Abstract}\stdspace\theabstract

{\bf AMS Classification}\stdspace\theprimaryclass
\ifx\thesecondaryclass\relax\else; \thesecondaryclass\fi\par
{\bf Keywords}\stdspace \thekeywords\par}\vglue 7truept

}   

\ifplaintex
\font\phead=cmsl9 scaled 950
\font\pnum=cmbx10 scaled 913
\font\pfoot=cmsl9 scaled 950
\headline{\vbox to 0pt{\vskip -4.5mm\line{\small\phead\ifnum
\count0=\startpage ISSN 1472-2739 (on-line) 1472-2747 (printed)
\hfill {\pnum\folio}\else\ifodd\count0\def\\{ }%
\ifx\theshorttitle\relax\thetitle\else\theshorttitle\fi\hfill{\pnum\folio}
\else\def\\{ and }{\pnum\folio}\hfill\ifx\theshortauthors\relax\theauthors
\else\theshortauthors\fi\fi\fi}\vss}}
\footline{\vbox to 0pt{\vglue 0mm\line{\small\pfoot\ifnum\count0=\startpage
\copyright\ \gtp\hfill\else
\agt, Volume \thevolumenumber\ (\thevolumeyear)\hfill\fi}\vss}}
\else
\font=cmsl9 scaled 1050
\font\lnum=cmbx10 
\font=cmsl9 scaled 1050
\makeatletter
\def\@oddhead{{\small\lhead\ifnum\count0=\startpage ISSN 1472-2739 
(on-line) 1472-2747 (printed)\hfill {\lnum\number\count0}\else\ifodd\count0
\def\\{ }\ifx\theshorttitle\relax \thetitle \else\theshorttitle\fi\hfill
{\lnum\number\count0}\else\def\\{ and }{\lnum\number\count0}
\hfill\ifx\theshortauthors\relax 
\theauthors\else\theshortauthors\fi\fi\fi}}\def\@evenhead{\@oddhead}
\def\@oddfoot{\small\lfoot\ifnum\count0=\startpage\copyright\ \gtp\hfill\else
\agt, Volume \thevolumenumber\ (\thevolumeyear)\hfill\fi}
\def\@evenfoot{\@oddfoot}
\makeatother
\fi
\let\maketitlepage\makeagttitle
\let\makeshorttitle\maketitlepage
\let\maketitle\maketitlepage

\newwrite\gtoutfile
\long\gdef\makeheadfile{  
{\def\\{, }\def\s{ }
\immediate\openout\gtoutfile head.xxx
\immediate\write\gtoutfile{To: math@arxiv.org}
\immediate\write\gtoutfile{Subject: put OR rep NNNNN:ppppp}
\immediate\write\gtoutfile{--text follows this line--}
\immediate\write\gtoutfile{Proxy-for: \ifx\theasciiauthors\relax
\theauthors\else\theasciiauthors\fi\s<\ifx\theasciiemail\relax\theemail\else\theasciiemail\fi>}
\immediate\write\gtoutfile{\noexpand\\}
\immediate\write\gtoutfile{Authors: \ifx\theasciiauthors\relax
\theauthors\else\theasciiauthors\fi}
{\def\\{ }\immediate\write\gtoutfile{Title: \ifx\theasciititle\relax
\thetitle\else\theasciititle\fi}}
\immediate\write\gtoutfile{Subj-class: GT or SG, GR etc}
\immediate\write\gtoutfile{MSC-class: \theprimaryclass\ifx\thesecondaryclass\relax\else, \thesecondaryclass\fi}
\immediate\write\gtoutfile{Journal-ref: Algebr. Geom. Topol. \thevolumenumber\s
(\thevolumeyear) \startpage-\finishpage}
\immediate\write\gtoutfile{Comments: Published by Algebraic and
Geometric Topology at}
\immediate\write\gtoutfile{\s\s\s  http://www.maths.warwick.ac.uk/agt/AGTVol\thevolumenumber/agt-\thevolumenumber-\thepapernumber.abs.html}
\immediate\write\gtoutfile{\noexpand\\}
\immediate\write\gtoutfile{}
\ifx\theasciiabstract\relax
\immediate\write\gtoutfile{\theabstract}\else
\immediate\write\gtoutfile{\theasciiabstract}\fi
\immediate\write\gtoutfile{}
\immediate\write\gtoutfile{\noexpand\\}
\immediate\write\gtoutfile{}
\immediate\closeout\gtoutfile}}  

\def\maketitlepage{\makeagttitle\makeheadfile}
\let\makeshorttitle\maketitlepage
\let\maketitle\maketitlepage

\def\ifplaintex{\expandafter\ifx\csname documentclass\endcsname\relax}

\def\gtp{{\mathsurround=0pt\it $\cal G\mskip-2mu$eometry \&\ 
$\cal T\!\!$opology $\cal P\!$ublications}}  

\def\recd{{\small Received:\qua\receiveddate\ifx\reviseddate\relax
\else\qquad Revised:\qua\reviseddate\fi\par}} 

\def\lognumber#1{\def\thelognumber{#1}}
\def\volumenumber#1{\def\thevolumenumber{#1}}
\def\volumeyear#1{\def\thevolumeyear{#1}}
\def\papernumber#1{\def\thepapernumber{#1}}
\def\pagenumbers#1#2{\def\startpage{#1}\def\finishpage{#2}}
\def\published#1{\def\publishdate{#1}}

\def\received#1{\def\receiveddate{#1}}
\def\revised#1{\def\reviseddate{#1}}
\def\accepted#1{\def\accepteddate{#1}}

\def\asciiaddress#1{\def\theasciiaddress{#1}}

\long\def\asciiabstract#1{\long\def\theasciiabstract{#1}}
\def\asciikeywords#1{\def\theasciikeywords{#1}}

\let\\\par\let\thelognumber\relax\let\thevolumenumber\relax
\let\thepapernumber\relax\let\thevolumeyear\relax\let\startpage\relax
\let\finishpage\relax\let\publishdate\relax\let\receiveddate\relax
\let\reviseddate\relax\let\accepteddate\relax\let\theasciititle\relax
\let\theasciiauthors\relax\let\theasciiaddress\relax
\let\theasciiabstract\relax\let\theasciikeywords\relax

\let\theasciiemail\relax

\ifplaintex
\font\logobig=cmssbx10 scaled 3836
\font\logomed=cmssbx10 scaled 2557
\else
\font\logobig=cmssbx10 scaled 4200
\font\logomed=cmssbx10 scaled 2800
\fi

\long\def\makeagttitle{   
\count0=\startpage
\agt\hfill      
\hbox to 45truept{\vbox to 0pt{\vglue -13truept{\logomed A\kern -.37em{\logobig 
T}\kern -.38em G}\vss}\hss}
\break
{\small Volume \thevolumenumber\ (\thevolumeyear)
\startpage--\finishpage\nl
Published: \publishdate}

\vglue .25truein

{\parskip=0pt\leftskip 0pt plus
1fil\def\\{\par\smallskip}{\Large\bf\thetitle}\par\medskip} \vglue
0.05truein

{\parskip=0pt\leftskip 0pt plus 1fil\def\\{\par}{\sc\theauthors}
\par\medskip}%
 
\vglue 0.03truein 

{\small\leftskip 25truept\rightskip 25truept{\bf Abstract}\stdspace\theabstract

{\bf AMS Classification}\stdspace\theprimaryclass
\ifx\thesecondaryclass\relax\else; \thesecondaryclass\fi\par
{\bf Keywords}\stdspace \thekeywords\par}\vglue 7truept

}   

\ifplaintex
\font\phead=cmsl9 scaled 950
\font\pnum=cmbx10 scaled 913
\font\pfoot=cmsl9 scaled 950
\headline{\vbox to 0pt{\vskip -4.5mm\line{\small\phead\ifnum
\count0=\startpage ISSN 1472-2739 (on-line) 1472-2747 (printed)
\hfill {\pnum\folio}\else\ifodd\count0\def\\{ }%
\ifx\theshorttitle\relax\thetitle\else\theshorttitle\fi\hfill{\pnum\folio}
\else\def\\{ and }{\pnum\folio}\hfill\ifx\theshortauthors\relax\theauthors
\else\theshortauthors\fi\fi\fi}\vss}}
\footline{\vbox to 0pt{\vglue 0mm\line{\small\pfoot\ifnum\count0=\startpage
\copyright\ \gtp\hfill\else
\agt, Volume \thevolumenumber\ (\thevolumeyear)\hfill\fi}\vss}}
\else
\font=cmsl9 scaled 1050
\font\lnum=cmbx10 
\font=cmsl9 scaled 1050
\makeatletter
\def\@oddhead{{\small\lhead\ifnum\count0=\startpage ISSN 1472-2739 
(on-line) 1472-2747 (printed)\hfill {\lnum\number\count0}\else\ifodd\count0
\def\\{ }\ifx\theshorttitle\relax \thetitle \else\theshorttitle\fi\hfill
{\lnum\number\count0}\else\def\\{ and }{\lnum\number\count0}
\hfill\ifx\theshortauthors\relax 
\theauthors\else\theshortauthors\fi\fi\fi}}\def\@evenhead{\@oddhead}
\def\@oddfoot{\small\lfoot\ifnum\count0=\startpage\copyright\ \gtp\hfill\else
\agt, Volume \thevolumenumber\ (\thevolumeyear)\hfill\fi}
\def\@evenfoot{\@oddfoot}
\makeatother
\fi
\let\maketitlepage\makeagttitle
\let\makeshorttitle\maketitlepage
\let\maketitle\maketitlepage

\newwrite\gtoutfile
\long\gdef\makeheadfile{  
{\def\\{, }\def\s{ }
\immediate\openout\gtoutfile head.xxx
\immediate\write\gtoutfile{To: math@arxiv.org}
\immediate\write\gtoutfile{Subject: put OR rep NNNNN:ppppp}
\immediate\write\gtoutfile{--text follows this line--}
\immediate\write\gtoutfile{Proxy-for: \ifx\theasciiauthors\relax
\theauthors\else\theasciiauthors\fi\s<\ifx\theasciiemail\relax\theemail\else\theasciiemail\fi>}
\immediate\write\gtoutfile{\noexpand\\}
\immediate\write\gtoutfile{Authors: \ifx\theasciiauthors\relax
\theauthors\else\theasciiauthors\fi}
{\def\\{ }\immediate\write\gtoutfile{Title: \ifx\theasciititle\relax
\thetitle\else\theasciititle\fi}}
\immediate\write\gtoutfile{Subj-class: GT or SG, GR etc}
\immediate\write\gtoutfile{MSC-class: \theprimaryclass\ifx\thesecondaryclass\relax\else, \thesecondaryclass\fi}
\immediate\write\gtoutfile{Journal-ref: Algebr. Geom. Topol. \thevolumenumber\s
(\thevolumeyear) \startpage-\finishpage}
\immediate\write\gtoutfile{Comments: Published by Algebraic and
Geometric Topology at}
\immediate\write\gtoutfile{\s\s\s  http://www.maths.warwick.ac.uk/agt/AGTVol\thevolumenumber/agt-\thevolumenumber-\thepapernumber.abs.html}
\immediate\write\gtoutfile{\noexpand\\}
\immediate\write\gtoutfile{}
\ifx\theasciiabstract\relax
\immediate\write\gtoutfile{\theabstract}\else
\immediate\write\gtoutfile{\theasciiabstract}\fi
\immediate\write\gtoutfile{}
\immediate\write\gtoutfile{\noexpand\\}
\immediate\write\gtoutfile{}
\immediate\closeout\gtoutfile}}  

\def\maketitlepage{\makeagttitle\makeheadfile}
\let\makeshorttitle\maketitlepage
\let\maketitle\maketitlepage

\def\ifplaintex{\expandafter\ifx\csname documentclass\endcsname\relax}

\def\gtp{{\mathsurround=0pt\it $\cal G\mskip-2mu$eometry \&\ 
$\cal T\!\!$opology $\cal P\!$ublications}}  

\def\recd{{\small Received:\qua\receiveddate\ifx\reviseddate\relax
\else\qquad Revised:\qua\reviseddate\fi\par}} 

\def\lognumber#1{\def\thelognumber{#1}}
\def\volumenumber#1{\def\thevolumenumber{#1}}
\def\volumeyear#1{\def\thevolumeyear{#1}}
\def\papernumber#1{\def\thepapernumber{#1}}
\def\pagenumbers#1#2{\def\startpage{#1}\def\finishpage{#2}}
\def\published#1{\def\publishdate{#1}}

\def\received#1{\def\receiveddate{#1}}
\def\revised#1{\def\reviseddate{#1}}
\def\accepted#1{\def\accepteddate{#1}}

\def\asciiaddress#1{\def\theasciiaddress{#1}}

\long\def\asciiabstract#1{\long\def\theasciiabstract{#1}}
\def\asciikeywords#1{\def\theasciikeywords{#1}}

\let\\\par\let\thelognumber\relax\let\thevolumenumber\relax
\let\thepapernumber\relax\let\thevolumeyear\relax\let\startpage\relax
\let\finishpage\relax\let\publishdate\relax\let\receiveddate\relax
\let\reviseddate\relax\let\accepteddate\relax\let\theasciititle\relax
\let\theasciiauthors\relax\let\theasciiaddress\relax
\let\theasciiabstract\relax\let\theasciikeywords\relax

\let\theasciiemail\relax

\ifplaintex
\font\logobig=cmssbx10 scaled 3836
\font\logomed=cmssbx10 scaled 2557
\else
\font\logobig=cmssbx10 scaled 4200
\font\logomed=cmssbx10 scaled 2800
\fi

\long\def\makeagttitle{   
\count0=\startpage
\agt\hfill      
\hbox to 45truept{\vbox to 0pt{\vglue -13truept{\logomed A\kern -.37em{\logobig 
T}\kern -.38em G}\vss}\hss}
\break
{\small Volume \thevolumenumber\ (\thevolumeyear)
\startpage--\finishpage\nl
Published: \publishdate}

\vglue .25truein

{\parskip=0pt\leftskip 0pt plus
1fil\def\\{\par\smallskip}{\Large\bf\thetitle}\par\medskip} \vglue
0.05truein

{\parskip=0pt\leftskip 0pt plus 1fil\def\\{\par}{\sc\theauthors}
\par\medskip}%
 
\vglue 0.03truein 

{\small\leftskip 25truept\rightskip 25truept{\bf Abstract}\stdspace\theabstract

{\bf AMS Classification}\stdspace\theprimaryclass
\ifx\thesecondaryclass\relax\else; \thesecondaryclass\fi\par
{\bf Keywords}\stdspace \thekeywords\par}\vglue 7truept

}   

\ifplaintex
\font\phead=cmsl9 scaled 950
\font\pnum=cmbx10 scaled 913
\font\pfoot=cmsl9 scaled 950
\headline{\vbox to 0pt{\vskip -4.5mm\line{\small\phead\ifnum
\count0=\startpage ISSN 1472-2739 (on-line) 1472-2747 (printed)
\hfill {\pnum\folio}\else\ifodd\count0\def\\{ }%
\ifx\theshorttitle\relax\thetitle\else\theshorttitle\fi\hfill{\pnum\folio}
\else\def\\{ and }{\pnum\folio}\hfill\ifx\theshortauthors\relax\theauthors
\else\theshortauthors\fi\fi\fi}\vss}}
\footline{\vbox to 0pt{\vglue 0mm\line{\small\pfoot\ifnum\count0=\startpage
\copyright\ \gtp\hfill\else
\agt, Volume \thevolumenumber\ (\thevolumeyear)\hfill\fi}\vss}}
\else
\font=cmsl9 scaled 1050
\font\lnum=cmbx10 
\font=cmsl9 scaled 1050
\makeatletter
\def\@oddhead{{\small\lhead\ifnum\count0=\startpage ISSN 1472-2739 
(on-line) 1472-2747 (printed)\hfill {\lnum\number\count0}\else\ifodd\count0
\def\\{ }\ifx\theshorttitle\relax \thetitle \else\theshorttitle\fi\hfill
{\lnum\number\count0}\else\def\\{ and }{\lnum\number\count0}
\hfill\ifx\theshortauthors\relax 
\theauthors\else\theshortauthors\fi\fi\fi}}\def\@evenhead{\@oddhead}
\def\@oddfoot{\small\lfoot\ifnum\count0=\startpage\copyright\ \gtp\hfill\else
\agt, Volume \thevolumenumber\ (\thevolumeyear)\hfill\fi}
\def\@evenfoot{\@oddfoot}
\makeatother
\fi
\let\maketitlepage\makeagttitle
\let\makeshorttitle\maketitlepage
\let\maketitle\maketitlepage

\newwrite\gtoutfile
\long\gdef\makeheadfile{  
{\def\\{, }\def\s{ }
\immediate\openout\gtoutfile head.xxx
\immediate\write\gtoutfile{To: math@arxiv.org}
\immediate\write\gtoutfile{Subject: put OR rep NNNNN:ppppp}
\immediate\write\gtoutfile{--text follows this line--}
\immediate\write\gtoutfile{Proxy-for: \ifx\theasciiauthors\relax
\theauthors\else\theasciiauthors\fi\s<\ifx\theasciiemail\relax\theemail\else\theasciiemail\fi>}
\immediate\write\gtoutfile{\noexpand\\}
\immediate\write\gtoutfile{Authors: \ifx\theasciiauthors\relax
\theauthors\else\theasciiauthors\fi}
{\def\\{ }\immediate\write\gtoutfile{Title: \ifx\theasciititle\relax
\thetitle\else\theasciititle\fi}}
\immediate\write\gtoutfile{Subj-class: GT or SG, GR etc}
\immediate\write\gtoutfile{MSC-class: \theprimaryclass\ifx\thesecondaryclass\relax\else, \thesecondaryclass\fi}
\immediate\write\gtoutfile{Journal-ref: Algebr. Geom. Topol. \thevolumenumber\s
(\thevolumeyear) \startpage-\finishpage}
\immediate\write\gtoutfile{Comments: Published by Algebraic and
Geometric Topology at}
\immediate\write\gtoutfile{\s\s\s  http://www.maths.warwick.ac.uk/agt/AGTVol\thevolumenumber/agt-\thevolumenumber-\thepapernumber.abs.html}
\immediate\write\gtoutfile{\noexpand\\}
\immediate\write\gtoutfile{}
\ifx\theasciiabstract\relax
\immediate\write\gtoutfile{\theabstract}\else
\immediate\write\gtoutfile{\theasciiabstract}\fi
\immediate\write\gtoutfile{}
\immediate\write\gtoutfile{\noexpand\\}
\immediate\write\gtoutfile{}
\immediate\closeout\gtoutfile}}  

\def\maketitlepage{\makeagttitle\makeheadfile}
\let\makeshorttitle\maketitlepage
\let\maketitle\maketitlepage